\documentclass[a4paper,11pt]{amsart}
\usepackage{amsmath,amssymb,amsthm,amscd}
\usepackage[]{latexsym,amssymb,amsmath,amsfonts, amsthm, verbatim}
\usepackage[all,cmtip,ps]{xy}
\usepackage[latin1]{inputenc}
\usepackage[english]{babel}

\usepackage{graphicx}
\DeclareMathOperator{\Img}{Im}
\newcommand{\Scal}{\ensuremath{\mathcal{S}}}
\newcommand{\Ucal}{\ensuremath{\mathcal{U}}}
\newcommand{\Xcal}{\ensuremath{\mathcal{X}}}
\newcommand{\perpe}[1]{\ensuremath{{\mathcal{X}}_{#1}}}

\let\Im\undefined
\let\End\undefined

\DeclareMathOperator{\Im}{Im}

\DeclareMathOperator{\Hom}{Hom}

\DeclareMathOperator{\End}{End}

\DeclareMathOperator{\Ker}{Ker}
\DeclareMathOperator{\Coker}{Coker}

\DeclareMathOperator{\Add}{Add}

\DeclareMathOperator{\Gen}{Gen}

\DeclareMathOperator{\Mod}{Mod}
\DeclareMathOperator{\modfg}{mod}

\def\Ext#1#2#3#4{\mathop{{\mathrm {Ext}}^{#1}_{#2}(#3,#4)}}
\def\Tor#1#2#3#4{\mathop{{\mathrm {Tor}}^{#1}_{#2}(#3,#4)}}

\def\pdim{\mathop{{\rm pd}}}

\def\dualita#1#2{\mathrel{
                 \mathop{\vcenter{
                 \offinterlineskip
                 \hbox to 1.2truecm{\rightarrowfill}
                 \hbox to 1.2truecm{\leftarrowfill}}}%
                 \limits_{#2}^{#1}}}

\renewcommand*{\epsilon}{\varepsilon}
\renewcommand*{\theta}{\vartheta}

\newtheorem{TEO}{Theorem}[section]
\newtheorem{PROP}[TEO]{Proposition}
\newtheorem{COR}[TEO]{Corollary}
\newtheorem{LEMMA}[TEO]{Lemma}
\theoremstyle{definition}
\newtheorem{DEF}[TEO]{Definition}
\newtheorem{ES}[TEO]{Example}
\newtheorem{REM}[TEO]{Remark}

\def\black{\ {\hbox{\vrule width 4pt height 4pt depth
0pt}}}
\def\fine{\ \black\vskip.4truecm}
\def\Proof{{\sl Proof.}\quad}

\def\co{\mathbin{\scriptscriptstyle\circ}}

\def\seb#1#2#3#4#5{0\rightarrow#1\stackrel{#4}{\rightarrow}#2
\stackrel{#5}{\rightarrow}#3\rightarrow 0}

\begin{document}

\title{ Tilting modules and universal localization}
\author{Lidia Angeleri H\" ugel, Maria Archetti}
\address{Lidia Angeleri  H\" ugel
\\ Dipartimento di Informatica - Settore Matematica
\\ Universit\`a degli Studi di Verona
\\ Strada Le Grazie 15 - Ca' Vignal 2
\\ 37134 Verona, Italy}
\email{lidia.angeleri@univr.it}
\address{Maria Archetti
\\ Department of Decision Sciences 
\\ Universit\`a Bocconi
\\ Via Roentgen 1
\\ 20136 Milan, Italy }
\email{maria.archetti@unibocconi.it}
\thanks{First named author acknowledges
support by  MIUR, PRIN-2008 "Anelli, algebre, moduli e categorie", by Progetto di Ateneo CPDA071244 of the University of Padova, and by  the DGI and the
European Regional Development Fund, jointly, through Project
 MTM2005--00934,  and by the Comissionat per Universitats i Recerca
of the Generalitat de Ca\-ta\-lunya, Project 2005SGR00206. }

\date{\today}

\begin{abstract}
We show that every tilting module of projective dimension one over a ring $R$  is associated in a natural way to the universal localization $R\to R_{\Ucal}$ at a set $\Ucal$
 of finitely presented modules of projective dimension one. We then investigate tilting modules of the form $R_\Ucal\oplus R_\Ucal/R$. Furthermore, we discuss the relationship between universal localization and the localization $R\to Q_{\mathcal G}$ given
 by a perfect Gabriel topology $\mathcal G$. Finally, we give some applications to Artin algebras and to Pr\"ufer domains.
 \end{abstract}

\maketitle

\section*{Introduction}
Tilting modules of projective dimension one are often constructed via a localization.
For example, if $\Sigma$ is a
{left Ore set of regular elements} in a ring $R$ with the property that the localization $\Sigma^{-1}R$ is an $R$-module of projective dimension at most one,
then $\Sigma^{-1}R\oplus \Sigma^{-1}R/R$ is a tilting right
$R$-module, see \cite{AHT, S, S2}.
More generally,  it was recently shown in  \cite{ASJ} that every injective homological ring epimorphism $R\to S$
such that $S_R$ has projective dimension at most one gives rise to a tilting $R$-module $S\oplus S/R$.

Note, however, that in general not all tilting modules arise as above from an
injective homological ring epimorphism.  For example,
 if  $R$ is  a commutative domain whose ring of fractions has projective
dimension at least two, then the {Fuchs' divisible module}
$\delta$ is a tilting $R$-module which is not of the form $S\oplus S/R$, cf.~Example~\ref{Fuchs}.

On the other hand, every tilting module $T$ of projective dimension one
is associated in a natural way to a reflective and coreflective subcategory of Mod$R$, which is obtained as
  perpendicular category
$$ {\mathcal X}_{T_1}=\{M_R\,|\ \Hom_R(T_1, M)=0=\Ext{1}{R}{T_1}{M}\}$$
of a certain module $T_1$ in the additive closure $\Add T$. By  a result of Gabriel and de la Pe\~na \cite[1.2]{GP}, the category  ${\mathcal X}_{T_1}$ is then associated to a ring epimorphism $\lambda:R\to S$. We show that, choosing  $T_1$ appropriately, one can find a set $\mathcal U$ of finitely presented modules of projective dimension one such that $\lambda$ is the universal localization $R\to R_{\mathcal U}$ of $R$ at $\mathcal U$ in the sense of Schofield \cite{Sc}. More precisely, we prove the following result.

\medskip

{\bf Theorem~\ref{teo:univlocalization1}} \emph{Let $T$ be a tilting module of projective dimension one. Then there are an exact sequence $0\to R\to T_0\to T_1\to 0$  and  a  set $\mathcal U$ of finitely presented modules of projective dimension one such that
\begin{enumerate}
\item  $T_0,T_1\in \Add T$,
\item  $\Gen T=\{M\in\Mod R\,|\, \Ext{1}{R}{U}{M}=0\ {\rm for\ every}\ U\in \mathcal U\}$,
\item
$ {\mathcal X}_{T_1}$ is equivalent to the category  of right $R_{\mathcal U}$-modules.
\end{enumerate}}

\smallskip

As a consequence, we see for instance that over an Artin algebra every finitely generated tilting module $T$  which is
of the form $S\oplus S/R$ for some injective homological ring epimorphism $\lambda:R\to S$ even arises from   universal localization $R\to R_{\mathcal U}$ at a set of finitely presented modules (Corollary~\ref{Artin}).

\medskip

We also study tilting modules arising from perfect localization. In particular, we describe which tilting modules  of the form $S\oplus S/R$ arise from the localization $R\to Q_{\mathcal G}$ given by a perfect Gabriel topology $\mathcal G$.

\smallskip

\textbf{Theorem \ref{teo:tiltingperfectloc}.}
\emph{Let $R$ be a ring and let $T_R$ be a tilting module of projective dimension one. The following
conditions are equivalent:
\begin{enumerate}
    \item There is a perfect Gabriel topology $\mathcal G$ such
    that $R$ embeds in $Q_{\mathcal G}$ and $Q_{\mathcal G}\oplus Q_{\mathcal G}/R$ is a tilting module
    equivalent to $T$.
    \item There is an exact sequence $\seb{R}{T_0}{T_1}{}{}$ such that $T_0,T_1\in\Add T$, $\Hom_R(T_1, T_0)=0$, and $\mathcal
    X_{T_1}$ is a Giraud subcategory of $\Mod R$.
\end{enumerate}
}

\smallskip

Observe that if $R$ is semihereditary, then every perfect localization $R\to Q_{\mathcal G}$ arises from universal localization at   a set of finitely presented modules (Proposition \ref{perfectisuniversal}). Over a Pr\"ufer domain, there is a converse result: every universal localization at a set of finitely presented cyclic modules can be viewed as the localization given by a perfect Gabriel topology (Proposition \ref{universalisperfect}).

\medskip

We  apply these results to investigate tilting modules over Pr\"ufer domains. Here the tilting classes are in one-one-correspondence with  perfect Gabriel topologies, as   shown by Bazzoni, Eklof and Trlifaj in \cite{BET}. More precisely, every tilting module $T$ is associated to  a
perfect Gabriel topology $\mathcal L$ such that the tilting class $\Gen T$ coincides with the class of $\mathcal L$-divisible modules. Moreover, if the localization $Q_{\mathcal L}$ has projective dimension at most one over $R$, then it was shown by Salce \cite{S2} that $T$ is equivalent to $ Q_{\mathcal L}\oplus Q_{\mathcal L}/R$. We recover Salce's result as a consequence of Theorem~\ref{teo:tiltingperfectloc}. Moreover,  we obtain that over a Pr\"ufer domain every tilting module
of the form $S\oplus S/R$   arises from  a universal localization $R\to R_{\mathcal U}$, as well as from a perfect localization $R\to Q_{\mathcal G}$, see Theorem~\ref{teo:14}.

\bigskip
{\bf Acknowldgements:} We thank Simion Breaz and Septimiu Crivei for drawing our attention on   perfect Gabriel topologies.

\section{Preliminaries}
\subsection*{I. Notation.} Let $R$ be a ring, and let $\Mod R$ be the category of all right $R$-modules. By a \emph{subcategory}  of $\Mod R$ we always mean a full subcategory which is closed under isomorphic images and direct summands.

We denote by $\modfg R$ the subcategory of modules
possessing a projective resolution consisting of finitely generated modules.

\medskip

 Given a class of modules $\mathcal C$, we denote
$$\mathcal C ^o = \{ M \in \Mod R \mid \Hom_R(C, M) = 0 \mbox{ for all } C \in
\mathcal C \},$$
$$\mathcal C^{\perp} = \{ M \in \Mod R \mid \Ext{i}{R}{C}{M} = 0 \mbox{ for all } C \in
\mathcal C \mbox{ and all } i > 0 \}.$$
The classes $^o \mathcal C$,
 and $^{\perp}\mathcal C$ are defined similarly.
The (right) \emph{perpendicular category} of $\mathcal C$ is denoted by
$$\mathcal X_{{\mathcal C}}=\mathcal C ^o\cap\mathcal C ^\perp$$

Finally, we denote by
$\Add \mathcal C $
the class consisting of all modules isomorphic to direct summands of
direct sums of modules of ${\mathcal C } $
and by $\Gen\mathcal C $
the class of modules generated
by modules of ${\mathcal C } $.

\bigskip

\subsection*{II.~Reflections} We start by recalling the notion of a reflective subcategory.

\begin{DEF}~\label{def:reflective}
Let $M$ be a right $R$-module and $\mathcal{C}$ a subcategory of $\Mod R$.

  A morphism
$f\in\Hom_R(M,C)$ with $C\in\mathcal{C}$ is said to be a
$\mathcal{C}$-\emph{preenvelope} of $M$ provided the morphism of
abelian groups
$\Hom_R(f,C')\colon\Hom_R(C,C')\rightarrow\Hom_R(M,C')$  is
surjective for each $C'\in\mathcal{C},$ that is, for each morphism
$f'\colon M\rightarrow C'$ there is a morphism $g\colon
C\rightarrow C'$ such that the following diagram is commutative.
\begin{eqnarray*}
\xymatrix@C=0.8cm@R=0.8cm{ M\ar[r]^f\ar[dr]_{f'} & C\ar@{-->}[d]^g
\\ & C'}
\end{eqnarray*}

Furthermore, a $\mathcal{C}$-preenvelope $f\in\Hom_R(M,C)$ is said to be a
$\mathcal{C}$-\emph{reflection} of $M$  provided the morphism of
abelian groups
$\Hom_R(f,C')\colon\Hom_R(C,C')\rightarrow\Hom_R(M,C')$ is
bijective for each $C'\in\mathcal{C},$ that is, the morphism
$g\colon C\rightarrow C'$ in the diagram above is always uniquely
determined.
In this case $f$ is also a
{$\mathcal{C}$-\emph{envelope}}, that is, a $\mathcal{C}$-preenvelope with the additional property that
  every $g\in\End_R(C)$ such that
$f=gf$  is an automorphism.

Finally,
 $\mathcal{C}$ is said to be a \emph{reflective} subcategory of $\Mod R$ if every $R$-module  admits a $\mathcal{C}$-reflection.
 \emph{Coreflective subcategories} are defined dually. A full subcategory
$\mathcal C$ of $\Mod R$ which is both reflective and coreflective is called \emph{bireflective}.
\end{DEF}

\begin{REM}~\label{rem:reflectivesub}
It is well known that a subcategory
$\mathcal{C}$ is a reflective subcategory of $\Mod R$ if and only
if the inclusion functor $\iota:
\mathcal C\hookrightarrow \Mod R$ has a left adjoint functor $\ell : \Mod R
\to \mathcal C$. In this case  a $\mathcal
C$-reflection of $M$ is given as $$\eta_M:M\to \iota\ell(M)$$ by the unit of the adjunction
$\eta:1_{\Mod R}\to \iota\,\ell$, see
 \cite[Chapter X, \S 1]{St}.
\end{REM}

Bireflective subcategories are closely related to ring epimorphisms.

\begin{DEF}
A ring homomorphism  $\lambda\colon
R\rightarrow S$ is called a \emph{ring epimorphism} if it is an epimorphism in the category of rings, that is, for every
pair of  morphisms of rings $\delta_i\colon S\rightarrow T,\
i=1,2,$ the condition $\delta_1\lambda=\delta_2\lambda$ implies
$\delta_1=\delta_2.$
Note that this holds true if and only if the restriction functor $$\lambda_\ast:\Mod S\to\Mod R$$ induced  by $\lambda$ is full, see \cite[Chapter~XI,
Proposition~1.2]{St}.

 Two ring epimorphisms
$\lambda :R\to S$ and $\lambda' :R\to S'$ are said to be \emph{equivalent} if there is
a ring isomorphism $\varphi :S\to S'$ such that $\lambda'=\varphi
\lambda$. The \emph{epiclasses} of $R$ are the equivalence classes
with respect to the equivalence relation defined above.
\end{DEF}

\begin{TEO}~\label{GP}\cite[1.2]{GP}, \cite{GL}, \cite[1.6.3]{I3}
The following assertions are
equivalent for
a subcategory $\mathcal X$   of $\Mod R$.
 \begin{enumerate}
\item $\mathcal X$ is a bireflective subcategory of $\Mod R$.
    \item $\mathcal X$ is is closed under isomorphic images, direct
sums, direct products, kernels and cokernels.
    \item There is a ring epimorphism $\lambda:R\to S$ such that $\mathcal X$ is the essential image of the restriction functor $\lambda_\ast:\Mod S\to\Mod R$.
\end{enumerate}
More precisely, there is a bijection between the epiclasses of the ring $R$ and the
bireflective subcategories of $\Mod R$.
Moreover, the map $\lambda:R\to S$ in condition (3), viewed as an $R$-homomorphism, is an $\mathcal X$-reflection of $R$.
\end{TEO}

 \medskip

 \subsection*{III.~Universal localization.}
Next, let us recall Schofield's notion of universal localization.

\begin{TEO}[{\cite[Theorem~4.1]{Sc}}]\label{def:universallocalization}
Let  $\Sigma$ be a set of morphisms between
finitely generated projective right $R$-modules. Then
there are a ring $R_\Sigma$ and a morphism of rings
$\lambda\colon R\rightarrow R_\Sigma$ such that
\begin{enumerate}
\item $\lambda$ is \emph{$\Sigma$-inverting,} i.e. if
$\alpha\colon P\rightarrow Q$ belongs to  $\Sigma$, then
$\alpha\otimes_R 1_{R_\Sigma}\colon P\otimes_R R_\Sigma\rightarrow
Q\otimes_R R_\Sigma$ is an isomorphism of right
$R_\Sigma$-modules, and
\item $\lambda$ is \emph{universal
$\Sigma$-inverting}, i.e. if $S$ is a ring such that there exists
a $\Sigma$-inverting morphism $\psi\colon R\rightarrow S$, then
there exists a unique morphism of rings $\bar{\psi}\colon
R_\Sigma\rightarrow S$ such that $\bar{\psi}\lambda=\psi$.
\end{enumerate}
\end{TEO}

The morphism $\lambda\colon R\rightarrow R_\Sigma$ is a ring epimorphism
with  $\Tor{R}{1}{R_\Sigma}{R_\Sigma}=0.$ It
 is called the \emph{universal localization of $R$ at
$\Sigma$}.

\medskip

Let now $\mathcal{U}$ be a set of  finitely presented right $R$-modules of projective dimension at most one.
For each $U\in\mathcal{U},$ consider a morphism $\alpha_U$ between
finitely generated projective right $R$-modules such that
$$0\to P\stackrel{{\alpha_U}}{\to} Q\to U\to 0$$
 We will denote by
$\lambda_{\Ucal}:R\to R_{\mathcal{U}}$ the universal localization of $R$ at the set
$\Sigma=\{\alpha_U\mid U\in\mathcal{U}\}.$ In fact,   $R_{\mathcal
U}$ does not depend on the class $\Sigma$  chosen,
cf.~\cite[Theorem~0.6.2]{Coh}, and we will also call it
the \emph{universal localization of $R$ at ${\mathcal{U}}$}.

\medskip

We now show that $\perpe{\Ucal}$ is the bireflective subcategory of $\Mod R$ corresponding to the ring epimorphism $\lambda_{\Ucal}$.

\begin{LEMMA}\label{RUmodules}
Let  $\mathcal{U}$ be a set of  finitely presented right $R$-modules of projective dimension at most one, and let $R_{\mathcal
U}$ be the universal localization of $R$ at ${\mathcal{U}}$.
Then every right $R_{\mathcal
U}$-module belongs to $\mathcal X_{\mathcal U}$.
\end{LEMMA}
\Proof
For each $U\in\mathcal{U},$ consider a sequence
$0\to P\stackrel{{\alpha_U}}{\to} Q\to U\to 0$ as above. Since $\alpha_U\otimes{}R_{\mathcal
U}$ is an isomorphism, we have $U\otimes{}R_{\mathcal
U}=\Tor{1}{R}{U}{R_{\mathcal
U}}=0$.

 Let now $M$ be a right $R_{\mathcal
U}$-module.  Note that the canonical map $\rho_M:M\to M\otimes_R{}R_{\mathcal
U}$ is an $R_{\mathcal
U}$-isomorphism, because  $\lambda\colon R\rightarrow R_{\mathcal
U}$ is a ring epimorphism.
For any $f\in\Hom_R(U, M)$ we then have a commutative diagram
$$\xymatrix{
U\ar[r]^{f}\ar[d]^{\rho _U}& M\ar[d]_{\rho _M}\\
{0=U\otimes _R{}R_{\mathcal U}}\ar[r]^-{f\otimes R_{\mathcal
U}}&M\otimes _R{}R_{\mathcal U}}$$ where $\rho _M{}f=0$ implies
$f=0$. Hence $M\in \mathcal U^0$.\\
Next, let $0\to M\stackrel{f}{\to} N\to U\to 0$ be an exact sequence in $\Mod R$. Then we have the
following commutative diagram with exact rows
$$\xymatrix{0\ar[r]&M\ar[r]^{f}\ar[d]^{\rho _M}&N\ar[r]\ar[d]_{\rho
_N}&U\ar[r]&0\\
0=\Tor{1}{R}{U}{R_{\mathcal U}}\ar[r]&M\otimes_R{}R_{\mathcal
U}\ar[r]_{f\otimes R_{\mathcal U}}&N\otimes_R{}R_{\mathcal
U}\ar[r]&U\otimes_R{}R_{\mathcal U}=0}$$ where $f\otimes R_{\mathcal
U}$ is an isomorphism. This implies that $\rho _N{}f$ is an isomorphism as well. Hence $f$ is a split monomorphism, and we deduce that $M\in \mathcal U^{\perp}$.
\fine

\smallskip

\begin{PROP}\label{teo:univlocalization2}
Let $\mathcal U$ be a set of of  finitely presented right $R$-modules of projective dimension at most one.
Then the following statements hold true.
\begin{enumerate}
\item The perpendicular category $\perpe{\Ucal}$ is bireflective.
\item $\perpe{\Ucal}$  coincides with the essential image of the restriction functor \linebreak
{$\Mod R_{\Ucal}\to\Mod R$}
induced by the  universal localization  at $\Ucal$.
\end{enumerate}
\end{PROP}
\Proof
(1) Clearly, $\perpe{\Ucal}$ is  closed  under direct products, and  $\mathcal U^0$ is closed under  direct products and submodules, hence also under direct sums. Furthermore,  the  assumptions on $\Ucal$ imply that $\mathcal
U^{\perp}$ is closed under epimorphic images and direct sums. So, we deduce  that  $\perpe{\Ucal}$ is closed under direct sums.

We now verify that  $\perpe{\Ucal}$ is closed under kernels.
Consider
$$\xymatrix{
0\ar[r]& \Ker f\ar[r]& Y\ar[rr]^{f}\ar@{->>}[dr]&& Z\\
&&& \Img f\ar@{^(->}[ur] }$$ with $Y, Z\in \perpe{\Ucal}$.
Since $\mathcal U ^{0}$ is closed under submodules and $\mathcal U
^{\perp}$ is closed under epimorphic images, we have $\Img f\in \mathcal U ^{0} \cap
\mathcal U ^{\perp}=  \perpe{\Ucal}$. Now, for $U\in \mathcal U$, applying
${\rm Hom}_R(U,{-})$ to the short exact
sequence $0\to{\Ker f}\to{Y}\to{\Img f}\to 0$, we get ${\rm Ext}^1_{R}({U},{{\Ker f}})=0$. This shows that  $\Ker f \in \perpe{\Ucal}$.

The closure  under cokernels is proved by similar arguments.

 So, we conclude from  Theorem~\ref{GP} that $\perpe{\Ucal}$ is bireflective.

(2)
 We know from Theorem \ref{GP} that there is a ring epimorphism
 {$\lambda:R\to S$} such that
 $\perpe{\Ucal}$ is the essential image of the restriction functor $\lambda_\ast:\Mod S\to\Mod R$
induced by $\lambda$. We claim that $\lambda$ is equivalent to the universal localization
 $\lambda_{\Ucal}:R\to R_{\Ucal}$  at $\Ucal$.
First of all, we choose a set $\Sigma=\{\alpha_U\mid U\in\mathcal{U}\}$ where the
$$0\to P_1\stackrel{{\alpha_U}}{\to} P_0\to U\to 0$$
 are exact sequences with finitely generated projective modules  $P_0,P_1$,
and we
claim that
 $\lambda$ is $\Sigma$-inverting.

Take $U\in\Ucal$ and set $\alpha=\alpha_U$. We have to show that $\alpha\otimes_R{}S$ is an isomorphism.
For any $S$-module $M$ we have $M_R\in \mathcal X_{\mathcal U}$, and thus we get the exact sequence
$$\xymatrix@-1.3pc{ { 0=\Hom_R(U, M)}\ar[r]&{\Hom_R({P_0}, M)}\ar[rrr]^{\Hom_R({\alpha}, M)}&&&{\Hom _R({P_1}, M)}\ar[r]&{\Ext{1}{R}{U}{M}=0}}$$
showing that
$\Hom_R(\alpha,  M)$ is an isomorphism.
Moreover, since $M\cong  \Hom_S(S, M)$ as $R$-modules, we have the following isomorphisms
$$0=\Hom_R(U, M) \cong \Hom_R(U, {\Hom_S(S, M)})\cong \Hom_S({U\otimes_R{}S}, M) .$$
In particular, choosing $M=U\otimes_R{}S$, we see
$\Hom_S({U\otimes_R{}S}, { U\otimes_R{}S})=0$, hence
$U\otimes_R{}S=\Coker (\alpha\otimes_R{}S)=0$.

Similarly, we
see that $${\Hom_S({\alpha\otimes_R{}S}, M}):\Hom_S({P_0\otimes_R{}S}, M)\to \Hom_S({P_1\otimes_R{}S}, M)$$
is an isomorphism.
In particular, choosing
$M=P_1\otimes_R{}S$, we obtain that $\Hom_S({\alpha\otimes_R{}S}, { P_1\otimes_R{}S})$ is an isomorphism and hence  $\alpha\otimes_R{}S$ is a split
monomorphism. Thus $\alpha\otimes_R{}S$ is an isomorphism.

Now, by the definition of  universal localization, there is a
(unique) map $\psi$ such that the following diagram commutes
$$\xymatrix{
R\ar[r]^{\lambda _{\mathcal U}}\ar[d]_{\lambda}&R_{\mathcal
U}\ar@{-->}[ld]^{\psi}\\
S}$$
Further, since
 $R_{\mathcal U}\in
\mathcal X_{\mathcal U}$ by Lemma \ref{RUmodules}, and $\lambda$ is an $\mathcal
X_{\mathcal U}$-reflection by Theorem~\ref{GP}, there is a (unique) map $\varphi$ such
that the following diagram commutes
$$\xymatrix{
R\ar[r]^{\lambda _{\mathcal U}}\ar[d]_{\lambda}&R_{\mathcal
U}\\
S\ar@{-->}[ru]_{\varphi}}$$
Now
$\psi\varphi\lambda=\psi\lambda_{\mathcal U}=\lambda$, hence $\psi\varphi=\mathrm{id}_R$. Moreover $\varphi\psi\lambda_{\mathcal
U}=\varphi\lambda=\lambda_{\mathcal U}$ and this implies
$\varphi\psi=\mathrm{id}_{R_{\mathcal U}}$. Hence we deduce that $\psi$ and $\varphi$
are isomorphisms, and the proof is complete.
\fine

\begin{REM}\cite{Pre,K1}
The map $\lambda:R\to R_{\mathcal U}$ can also be described as
 \emph{ring of definable scalars} for $\tilde{\Sigma}=\{P_1\oplus U\stackrel{(\alpha, 0)}{\to}P_0| \ \alpha\in \Sigma,\ U=\Coker
    \alpha\}$, or as \emph{biendomorphism ring} of a module $M$ which  is
    constructed as follows: take $N$ as the direct product of a
    representative set of the indecomposable pure-injective modules in $\mathcal X_{\mathcal
    U}$, set $\kappa=\mathrm{card} N$, and $M=N^{\kappa}$.
For details, see \cite[12.13, 12.16 and 11.7]{K1}.
\end{REM}

\bigskip

\subsection*{IV.~Tilting modules} Finally, let us review the notion of a tilting module and its relationship with ring epimorphisms.

\begin{DEF} A module $T$ is said to be  a {\em tilting module (of projective dimension at most one)}  if Gen$T=T^\perp$, or equivalently, if  the following conditions are
 satisfied:

  (T1) proj.dim$(T) \le 1$;

  (T2) ${\rm Ext}^1_R(T,{T^{(I)}}) = 0$ for each set $I$; and

  (T3) there is an    exact sequence
$0 \to R \to T_0 \to T_1 \to 0$ where  $T_0,T_1$ belong to Add$T$.

\smallskip

The class $T^\perp$ is then called a {\em tilting class}.
We  say that two tilting modules $T$ and $T'$ are {\em equivalent} if their tilting classes coincide.
\end{DEF}

\medskip

Here is a typical pattern for constructing tilting modules.

\begin{PROP}\cite[2.5]{ASJ}\label{prop:6}
 Let $\lambda:R\to S$ be an injective ring epimorphism with $\Tor{R}{1}{S}{S}=0$. Then $\pdim S_R\le 1$ if and only if $S\oplus S/R$ is a tilting right $R$-module.
\end{PROP}

 The following statements,
  relying on Theorem \ref{GP} and results from  \cite{CTT},  are shown in  \cite[proof of Theorem 2.10]{ASJ}.

\begin{LEMMA}\cite{ASJ}\label{reflective}
Let $T$ be a  tilting module  of projective dimension one, and let
 $0\to R \to T_0\to T_1\to 0$
 be an exact sequence with $T_0,T_1\in \Add T$. Then
 \begin{enumerate}
 \item $\Gen T=\Gen T_0=T_1\,^\perp$.
 \item $T_0\oplus T_1$ is a tilting module equivalent to $T$.
 \item $\mathcal X_{T_1}$ is a bireflective subcategory of $\Mod R$, so there is a ring epimorphism
 $\lambda:R\to S$ such that  $\mathcal X_{T_1}$ coincides with
the essential image of the restriction functor $\lambda_\ast:\Mod S\to\Mod R$
induced by $\lambda$.
 \end{enumerate}
 \end{LEMMA}

 In fact, the observations above are used to prove  the following result.

\begin{TEO}\cite[2.10]{ASJ}~\label{teo:tiltinglocalization}
Let $T_R$ be a tilting module  of projective dimension one. The following assertions are
equivalent:
\begin{itemize}
    \item [(1)] There is an injective ring epimorphism $\lambda :R\to S$ such that $\Tor{R}{1}{S}{S}=0$ and $S\oplus S/R$ is a tilting module equivalent to $T_R$.
    \item [(2)] There is an exact sequence $0\to{R}\stackrel{a}{\to}{T_0}\to{T_1}\to 0$ with $T_0, T_1\in \Add T$ and $\Hom_R(T_1, T_0)=0.$
\end{itemize}
Moreover, under these conditions, $a:R\to T_0$ is a $T^\perp$-envelope of $R$, and  $\mathcal X_{T_1}$ coincides with
the essential image of the restriction functor $\lambda_\ast:\Mod S\to\Mod R$
induced by $\lambda$.
\end{TEO}

 \bigskip

\section{Tilting modules arising from universal localization}

Aim of this section is to show that every tilting module of projective dimension one is associated in a natural way to a ring epimorphism which, moreover, can be interpreted as a universal localization at a set of finitely presented modules of projective dimension one.

\medskip

The following result by Bazzoni and Herbera will play an important role.
\begin{TEO} \cite{BH} \label{BH}
Let $T$ be a tilting  module of projective dimension one. Then there is a set $\mathcal S$  of modules in $\modfg R$ of projective dimension one such that
$T^\perp=\mathcal S^\perp$. More precisely, $\mathcal S$ can be chosen as a set of representatives of the isomorphism classes of non-projective modules from
${}^\perp(T^\perp)\cap\modfg R$.
\end{TEO}

Recall that, given a module $M$, an increasing chain of submodules
$\mathcal M = ( M_\alpha \mid \alpha \leq \sigma )$ of $M$, indexed by an ordinal $\sigma$,
 is called a {\em filtration} of $M$ provided that $M_0 = 0$,
$M_\alpha = \bigcup_{\beta < \alpha} M_\beta$ for all limit
ordinals $\alpha \leq \sigma$, and $M_\sigma = M$.
Moreover, if all consecutive factors $M_{\alpha + 1}/M_\alpha$, $\alpha < \sigma$, belong to a given subcategory
$\mathcal C$ of $\Mod R$,  we say that $M$ is {\em $\mathcal C$--filtered}.

\medskip

Let us now fix a tilting module $T$ of projective dimension one, and  let   $\Scal$ be a set of representatives of the isomorphism classes of the non-projective modules in  ${}^\perp(T^\perp)\cap\modfg R$. Then $T^\perp=\Scal^\perp$ by   Theorem \ref{BH}.
Hence by
\cite[3.2.1]{GT} there exists an
exact sequence
$0\to R\to T_0\to T_1\to 0$ where $T_0\in T^\perp$ and $T_1$
is $\Scal$-filtered.

\begin{TEO}\label{teo:univlocalization1}
There exist
an exact sequence $$0\to R \to T_0\to T_1\to 0$$
and a set $\mathcal U\subseteq\modfg R$ of modules of projective dimension one
such that
\begin{enumerate}
\item
$T_0,T_1\in\Add T$ and $T_1$ is $\Ucal$-filtered.
\item
$\Gen T=\Ucal^\perp$.
\item
$\perpe{T_1}=\perpe{\Ucal}$ coincides with the essential image of the restriction functor $\Mod R_{\Ucal}\to\Mod R$
induced by the  universal localization at $\Ucal$.
\end{enumerate}
\end{TEO}
\Proof
$(1)$ From the discussion above we know that there is an
exact sequence
$0\to R\to T_0\to T_1\to 0$ where $T_0\in T^\perp$ and $T_1$
is $\Scal$-filtered.

  The module $T_1$ then belongs to
${}^\perp(T^\perp)$ by \cite[3.1.2]{GT}. Since $ T^\perp=\Gen T$ is closed under quotients,
$T_1$ also belongs to $ T^\perp$. So $T_1\in  T^\perp\cap {}^\perp(T^\perp)=\Add T$.

Moreover, $T_0\in T^\perp\cap {}^\perp(T^\perp)=\Add T$, because $R$ belongs to $^\perp(T^\perp)$ which is closed under extensions.

Take now  an $\Scal$-filtration $( M_\alpha \mid \alpha \leq \sigma )$  of $T_1$ and
set  $$\Ucal=\{M_{\alpha+1}/M_\alpha\,\mid\, \alpha<\sigma\}.$$  Note that  $\Ucal$ consists of modules of projective dimension one by \cite[5.1.8]{GT}, and $T_1$ is obviously $\Ucal$-filtered.

\medskip

$(2)$ From
$\mathcal U\subseteq {}^\perp(T^\perp)$ we infer $\Gen T=T^\perp\subseteq \mathcal U^{\perp}$ .

For the reverse inclusion, recall from Lemma \ref{reflective} that $\Gen T=T_1\,^\perp$.
Since $T_1$ is $\mathcal U$-filtered and
$\mathcal U\subseteq{} ^{\perp}(\mathcal U^{\perp})$, we deduce from \cite[3.1.2]{GT} that $T_1\in{}^{\perp}(\mathcal U^{\perp})$, hence
$\Ucal^\perp\subseteq T_1\,^\perp$.

\medskip

$(3)$ We start by showing
$\mathcal X_{T_1}\subseteq \mathcal X_{\mathcal U}$. Let $X\in
\mathcal X_{T_1}$. Then $X\in {T_1}^{\perp}=\Gen T= \mathcal U^{\perp}$.
Assume $X\notin \mathcal
U^{0}$. Then there exists $0\ne f:U\to X$ for some $U=M_{\alpha
+1}/M_{\alpha}\in \mathcal U$. This implies that also $g:M_{\alpha
+1}\twoheadrightarrow U\to X$ is different from zero.

Indeed, for
all $\beta >\alpha$ there exists $0\ne g_{\beta}: M_{\beta}\to X$, as we are going to show.
For $\beta =\alpha +1$, we take $g_{\beta}=g$.
Given $g_{\beta}$, we consider
$$\xymatrix{
0\ar[r]&M_{\beta}\ar[rr]\ar[dr]_{g_{\beta}}&&M_{\beta
+1}\ar[r]\ar@{-->}[dl]^{g_{\beta +1}}&M_{\beta +1}/M_{\beta}\ar[r]&0\\
&&X}$$ and we use that the map $g_{\beta}$ extends to $g_{\beta
+1}$ since $X\in \mathcal U^{\perp}$.
Further, for  a limit ordinal $\beta$, we have that the $g_{\gamma}:M_{\gamma}\to X$
with $\gamma <\beta$ form a direct system inducing a non-zero map
$g_{\beta}:M_{\beta}=\cup_{\gamma<\beta}M_{\gamma}\to X$.

In particular, we obtain
$\Hom_R({M_{\sigma}}, { X})\ne 0$. But $M_{\sigma}=T_1$,  so
$X\notin
T_1^{0}$, a contradiction.
Thus we
conclude that $X\in \mathcal U^{\perp}\cap \mathcal U^{0}=\mathcal
X_{\mathcal U}$.

We now show
$\mathcal X_{\mathcal U}\subseteq \mathcal X_{T_1}$. Let
$X\in \mathcal X_{\mathcal U}=\mathcal U^{\perp}\cap \mathcal
U^0$. We already know that $X$ then belongs to $\Gen T=T_1\,^\perp$, so it remains to verify that $X\in T_1^0$. Since $T_1=M_{\sigma}$, this will follow once we show
 that
$\Hom_R({M_{\beta}}, {X})=0$ for all $\beta\le \sigma$.\\
The claim is clear for $\beta=0$ since $M_0=0$, and for $\beta =1$ since
$M_1\in \mathcal U.$ If $\beta=\alpha +1$ for some $\alpha$,  then we have an exact
sequence
$0\to {M_{\alpha}}\to{M_{\alpha+1}}\to{M_{\alpha+1}/M_{\alpha}}\to 0$ where $M_{\alpha+1}/M_{\alpha}\in \mathcal U\subseteq{} ^0X$, and   $M_{\alpha}$ belongs to $^0X$ by inductive assumption.  Since the
class $^0X$ is closed under extensions, we infer that also $M_{\alpha +1}$
belongs to $^0X$.\\
Finally, if $\beta$ is a limit ordinal, then
$M_{\beta}=\varinjlim_{\alpha<\beta}M_{\alpha}$, and again
by inductive assumption
$\Hom_R({M_{\beta}}, X) \cong\varprojlim \Hom_R({M_{\alpha}}, {X})=0$.

So the first claim is verified, and Proposition \ref{teo:univlocalization2} completes the proof.
\fine

\begin{DEF}
We will say that a tilting module $T$ \emph{arises from universal localization} if there is a set $\mathcal U\subset\modfg R$ of modules of projective dimension at most one such that  $R$ embeds in $R_{\mathcal U}$ and
$R_{\mathcal U}\oplus R_{\mathcal U}/R$ is a tilting $R$-module equivalent to $T$.
\end{DEF}

\begin{COR}~\label{cor:univlocalization}
Let $0\to R \to T_0\to T_1\to 0$ and $\mathcal U$ be as in
Theorem~\ref{teo:univlocalization1}. If $\lambda _{\mathcal U}$ is a
$\mathcal U^{\perp}$-preenvelope of $R$, then $T$ arises from universal localization.
\end{COR}
\Proof Since $\mathcal U^{\perp}$ contains all injective modules,  $\lambda _{\mathcal U}$ is injective.
Moreover,
if $\lambda_{\mathcal U}$ is a $\mathcal U^{\perp}$-preenvelope, then $\Gen R_{\mathcal
U}=\mathcal U^{\perp}=\Gen T$, cf.~\cite[3.12]{ASJ}.

  We claim that  pd$R_{\mathcal U}\le 1$. By assumption, there are $R$-epimorphisms
$$f: T^{(I)}\twoheadrightarrow R_{\mathcal U}\, \qquad \text{and} \qquad g: R_{\mathcal U}\,^{(J)}
\twoheadrightarrow T^{(I)}.$$ The composition $fg:
R_{\mathcal U}\,^{(J)}\twoheadrightarrow R_{\mathcal U}\,$ is an $R$-epimorphism, and also
an $R_{\mathcal U}\,$-epimorphism since $\Mod R_{\mathcal U}\,$ is a
full subcategory of $\Mod R$. Therefore $fg$ is a split $R_{\mathcal U}\,$-epimorphism, so there is $h\in \Hom_{R_{\mathcal U}}\,(R_{\mathcal U}\,, R_{\mathcal U}\,^{(J)})= \Hom_R(R_{\mathcal U}\,,
R_{\mathcal U}\,^{(J)})$ such that $fgh=id_{R_{\mathcal U}\,}$. Thus $f$
is a split $R$-epimorphism, showing that $\pdim R_{\mathcal U}\,\le \pdim T^{(I)}\le 1$.

 Since $\Tor{R}{1}{R_{\mathcal U}}{R_{\mathcal U}}=0$, we conclude  from
Proposition~\ref{prop:6} that $R_{\mathcal U}\oplus R_{\mathcal U}/R$ is a tilting
module equivalent to $T$. \fine

\begin{COR}~\label{cor:univlocfiltration}
  Assume that there is an exact sequence $\seb{R}{T_0}{T_1}{}{}$ such that  $T_0,T_1\in \Add T$,  $T_1$ is $\mathcal S$-filtered, and  $\Hom_R(T_1, T_0)=0$. Then $T$ arises from  universal localization.
\end{COR}
\Proof By Theorem~\ref{teo:tiltinglocalization} there is an injective ring epimorphism $\lambda:R\to S$ such that the $R$-module $S\oplus S/R$ is a tilting module equivalent to $T$, and moreover, $\mathcal X_{T_1}$ coincides with
the essential image of the restriction functor $\lambda_\ast:\Mod S\to\Mod R$
induced by $\lambda$. On the other hand, we have seen in Theorem~\ref{teo:univlocalization1} that $\perpe{T_1}$ coincides with the essential image of the restriction functor $\Mod R_{\Ucal}\to\Mod R$
induced by the  universal localization at a set $\Ucal$ of finitely presented modules of projective dimension one. Then
 it follows from   Theorem~\ref{GP} that
$\lambda:R\to S$ and $\lambda_{\mathcal U}:R\to R_{\mathcal U}$ are in the same epiclass. So, $\lambda_{\mathcal U}$ is injective and $R_{\mathcal U}\oplus R_{\mathcal U}/R$ is a tilting module equivalent to $T$. \fine

\begin{ES}~\label{es:pruferartin}
    Assume that  $T\in\modfg R$ and that the category    $T^{\perp}\cap\modfg R$ is covariantly finite in $\modfg R$, that is, every module in $\modfg R$ has a $T^{\perp}\cap\modfg R$-preenvelope. Assume further that there is an exact sequence $\seb{R}{T_0}{T_1}{a}{}$ such that $T_0,T_1\in \Add T$, and $\Hom_R(T_1, T_0)=0$. Then $T$ arises from universal localization.\\
        In fact, this will follow immediately from Corollary~\ref{cor:univlocfiltration} once we prove that $T_1$ belongs to $\modfg R$ (and is therefore trivially $\mathcal S$-filtered).\\
        Let
     us start by considering a $T^{\perp}\cap\modfg R$-preenvelope $f:R\to B$. We claim that $f$ is even a $T^{\perp}$-preenvelope. Indeed, if $h:R\to X$ with $X\in T^{\perp}=\Gen T$, then there exists an epimorphism $g:T^{(\alpha)}\twoheadrightarrow X$, and $h$ factors through $g$ via a homomorphism  $h':R\to T^{(\alpha)}$. Since the image of $h'$ is contained in a finite subsum $T^{(\alpha_0)}$ of $T^{(\alpha)}$, we can even factor  $h=g'\,h''$ where   $h'':R\to T^{(\alpha_0)}$ and $g:T^{(\alpha_0)}\twoheadrightarrow X$. Now $T^{(\alpha_0)}\in T^{\perp}\cap\modfg R$, so there is a map $\tilde{h}:B\to T^{(\alpha_0)}$ such that $h''=\tilde{h}f$, hence $h=g'\tilde{h}f$. This proves our claim.\\
      On the other hand, we know from Theorem~\ref{teo:tiltinglocalization} that $a$ is a $T^{\perp}$-envelope. Thus $T_0$ is isomorphic to a direct summand of $B$, and since $B\in\modfg R$, we infer that $T_0, T_1$ belong to $\modfg R$.
\end{ES}

 In particular, we deduce the following result  from \cite{AR}.
\begin{COR}\label{Artin}
Let $R$ be an  Artin algebra, and let $T$ be  a finitely generated tilting right $R$-module of projective dimension one. The following assertions are
equivalent:
\begin{itemize}
    \item [(1)] There is a set of finitely generated modules $\Ucal$ of projective dimension one such that  $R_\Ucal\oplus R_\Ucal/R$ is a tilting module equivalent to $T_R$.
    \item [(2)] There is an exact sequence $0\to{R}\stackrel{a}{\to}{T_0}\to{T_1}\to 0$ with $T_0, T_1\in \Add T$ and $\Hom_R(T_1, T_0)=0.$
\end{itemize}
\end{COR}

In the last section, we will see that a similar result holds true over Pr\"ufer domains.

\bigskip

\section{Tilting modules arising from perfect localization}

In this section we investigate  tilting modules arising from
perfect localization. We start by
 recalling some basic notions and results. For details we refer to \cite{G,P,St}.

\begin{DEF}~\label{def:giraud}
(1) A (full) subcategory $\mathcal X$ of $\Mod R$ is called a \emph{Giraud subcategory} if the canonical inclusion
$\iota : \mathcal X\to \Mod R$ has a left adjoint $\ell :\Mod R\to \mathcal X$ which is an exact functor.
The
composition functor $L=\iota \co \ell $ is then called
\emph{localization functor}.

\medskip

(2)
A non-empty set of right $R$-ideals $\mathcal G$ is said to be a \emph{Gabriel
topology} on $R$ if satisfies the following conditions:
\begin{enumerate}
    \item[(a)] If $I\in \mathcal G$ and a right ideal $K$ contains $I$, then $K$ belongs to $\mathcal G$.
    \item[(b)] If $I$ and $K$ belong to $\mathcal G$ then also $I\cap K$ belongs to $\mathcal G$.
    \item[(c)] If $I\in \mathcal G$ and $x\in R$ then $(I:x)=\{r\in R| \ xr\in I\}$ belongs to $\mathcal G$.
    \item[(d)] If $K$ is a right ideal and if there is some $I\in \mathcal G$ such that $(K:x)\in \mathcal G$ for any $x\in I$, then also $K$ belongs to $\mathcal G$.
\end{enumerate}
Further, a Gabriel topology $\mathcal G$ is \emph{of finite type} if it has a \emph{basis} of finitely generated ideals, that is, every $I\in\mathcal G$ contains a finitely generated right ideal $I'\in\mathcal G$.

\medskip

 (3)
Let  $\mathcal G$ be a Gabriel topology on $R$. A right $R$-module $C$ is said to be
$\mathcal G$-\emph{closed} if for any
short exact sequence
$\seb{I}{R}{R/I}{}{}$ {with} $I\in \mathcal G$
the morphism of abelian groups $\Hom_{R}(R, C)\to \Hom_{R}(I, C)$
is bijective.
\\
 A left $R$-module $_RX$ is said to be $\mathcal
    G$-\emph{divisible} if $IX=X$ for all $I\in \mathcal G$.

\medskip

 (4) A pair of subcategories $(\mathcal T, \mathcal F)$ is said to be a \emph{torsion pair} if ${\mathcal T}={}^o{\mathcal F}$ and  ${\mathcal T}^o={\mathcal F}$. In this case, $\mathcal T$ is a \emph{torsion class}, that is, it is closed under epimorphic images, extensions, and direct sums. If, in addition, $\mathcal T$ is closed under submodules, then $(\mathcal T, \mathcal F)$ is called  a \emph{hereditary} torsion pair.
\end{DEF}

\begin{TEO} \cite[VI, 5.1 and X, 2.1]{St}~\label{teo:hereditarygiraud}
There are
bijective correspondences between
the hereditary torsion pairs in $\Mod R$,  the Gabriel topologies on $R$, and the Giraud subcategories of $\Mod R$.

More precisely, under these bijections, a hereditary
torsion pair $(\mathcal T, \mathcal F)$ in $\Mod R$ corresponds
to the Gabriel topology  $$\mathcal G=\{ I\le R \,\mid\, R/I\in\mathcal T\}$$ as well as
to the Giraud subcategory $\mathcal X_{\mathcal T}$. Conversely, a Giraud subcategory $\mathcal X$ with localization functor $L$ is associated  to the hereditary torsion pair with torsion class
$$\mathcal T:=\{M\in \Mod R| \ L(M)=0\}$$
Finally, if  $\mathcal G$ is a Gabriel topology, then the category $\mathcal X({\mathcal G})$ of all $\mathcal G$-closed modules is the corresponding Giraud subcategory.\end{TEO}

\medskip

Let now $\mathcal G$ be a Gabriel topology. Consider  the adjoint pair $(\ell, \iota)$   corresponding to the Giraud subcategory $\mathcal X({\mathcal G})$ of all $\mathcal G$-closed modules, and the localization functor $L=\iota\co\ell$. Recall from Remark \ref{rem:reflectivesub} that the unit of the adjunction
$\eta_{M}: M\to L(M)$ defines an $\mathcal X({\mathcal G})$-reflection.
In particular,  $\eta_{R}: R\to L(R)$ induces a ring structure on  $Q_{\mathcal G}=L(R)$, and we obtain a ring homomorphism  $\lambda_{\mathcal G}:R\to Q_{\mathcal G}$.

\begin{TEO}[{\cite[XI, 3.4]{St}}]~\label{teo:perfectloc}
Let $\mathcal G$ be a Gabriel topology on $R$, and let  $\mathcal X({\mathcal G})$ be the corresponding  Giraud subcategory of all $\mathcal G$-closed modules. The following
assertions are equivalent.
\begin{enumerate}
    \item  $\mathcal X({\mathcal G})$  is a coreflective subcategory of $\Mod R$.
    \item $\mathcal X({\mathcal G})$  coincides with the essential image of the restriction functor $\Mod Q_{\mathcal G}\to\Mod R$ induced by $\lambda_{\mathcal G}$.
    \item The left $R$-module $Q_{\mathcal G}$ is $\mathcal G$-divisible.
\end{enumerate}
\end{TEO}

\begin{DEF} A Gabriel topology $\mathcal G$ that satisfies the  equivalent conditions of Theorem \ref{teo:perfectloc} is called a \emph{perfect Gabriel topology}.\end{DEF}

\begin{REM}~\label{finitetype}
Let  $\mathcal G$ be a  Gabriel topology on $R$, and let $(\mathcal T,\mathcal F)$ be the corresponding hereditary torsion pair.
 The torsion class  $\mathcal
T$ consists of all modules $X$ such that every $x\in
X$ has annihilator $\text{ann}_R(x)\in \mathcal G$, and the torsion-free class $\mathcal F$  is given by the modules $M$ for which the ${\mathcal X(G)}$-reflection $\eta_M:M\to L(M)$ is injective  \cite[VI, 5.1, and X, 1.5]{St}.

Assume now that
  $\mathcal G$ is perfect. Then  $\lambda_{\mathcal G}:R\to Q_{\mathcal G}$ is a ring epimorphism,   $Q_{ \mathcal G}$ is a flat left $R$-module, and $\mathcal G$  is  a Gabriel topology of finite type {\cite[XI, 3.4]{St}}. Moreover, $\mathcal F={\mathcal C}^o$ {where} $\mathcal C$ is the class of all finitely presented cyclic modules in $\mathcal T$.
This follows easily from {\cite[VI, 3.6]{St}} by using that $\mathcal G$ has finite type.
\end{REM}

We now fix a tilting module $T$ of projective dimension one
together with an exact sequence  $$\seb{R}{T_0}{T_1}{}{}$$ where $T_i\in \Add T$. We know from Lemma~\ref{reflective} that
$\mathcal X_{T_1}$ is a bireflective subcategory of $\Mod R$. We
denote by $\ell $ the left adjoint of the inclusion functor
$\iota:\mathcal X_{T_1}\hookrightarrow \Mod R$.
In \cite{CTT}, the  functor $\ell$ is constructed explicitly by using {Bongartz preenvelopes}.
More precisely, if $M_R$ is a right $R$-module, and $c$ is the minimal number of
generators of $\Ext{1}{R}{T_1}{M}$ as a module over $\End_R(T_1)$, then there exists an exact sequence
$$\seb{M}{M_0}{T_1^{(c)}}{i}{}$$
with $M_0\in T_1^{\perp}$. In particular, $i$ is a special $\Gen T$-preenvelope of $M$, called the \emph{Bongartz preenvelope} of $M$. It is shown in {\cite[1.3]{CTT}}
that $\ell(M)$ can be computed as
$$\ell (M):=M_0/\text{tr}_{T_1}M_0$$
where $\text{tr}_{T_1}M_0=\Sigma\{\Im f| \ f\in \Hom_R(T_1, M_0)\}$ denotes the $T_1$-trace of $M_0$.
We use this description in order to determine the kernel of the functor $\ell$.

\begin{LEMMA}~\label{lemma:7}
 For each $M\in \Mod R$ fix a Bongartz preenvelope $M_0$.  Then the following statements are equivalent.
 \begin{enumerate}
 \item $M\in {}^0\mathcal X_{T_1}$
 \item $\ell(M)=0$
 \item $ M_0\in \Gen T_1$
 \end{enumerate}
  Moreover, these conditions are satisfied whenever $M\in \Gen T_1$.
\end{LEMMA}
\Proof (1)$\Rightarrow$ (2): Recall that $\eta_M: M\to \iota\ell(M)$ is the $\mathcal X_{T_1}$-reflection of $M$, which is uniquely determined up to isomorphism. So, if $M\in {}^0\mathcal X_{T_1}$, we must have $\ell (M)=0$.

(2)$\Rightarrow$ (3) is clear.

(3)$\Rightarrow$ (1): Assume that $M_0\in \Gen T_1$ and $M\notin{} ^0\mathcal
    X_{T_1}$. Then there exists a map $0\ne f\in \Hom_R(M, X)$ for some
    $X\in \mathcal X_{T_1}\subseteq T_1^{\perp}$. Then
    $$\xymatrix{
    0\ar[r]&M\ar[r]^{i}\ar[d]_{f}&M_0\ar[r]\ar@{-->}[ld]^{h}&T_1^{(c)}\ar[r]&0\\
    &X}$$
    there is a map $0\ne h\in \Hom_R(M_0, X)$. Since $M_0\in \Gen
    T_1$ there exists a map $0\ne h'\in \Hom_R(T_1, X)$. But $X\in
    \mathcal X_{T_1}\subseteq T_1^0$, a contradiction.

    Finally,  if $M\in \Gen T_1\subseteq \Gen T=T_1^{\perp}$, then
    $\Ext{1}{R}{T_1^{(c)}}{M}=0$, so $M_0\cong M\oplus T_1^{(c)}\in \Gen
    T_1$.
    \fine

\begin{PROP}\label{prop:8}
Let $\lambda:R\to S$ be a  ring
epimorphism such that  the essential image of the restriction functor $\lambda_\ast:\Mod S\to \Mod R$ coincides with
 $\mathcal X_{T_1}$.  Then the following assertions are equivalent:
\begin{enumerate}
    \item $_RS$ is a flat left $R$-module.
    \item $\mathcal X_{T_1}$ is a Giraud subcategory of $\Mod R$.
    \item All submodules of modules in $\Gen T_1$ belong to  $^0\mathcal X_{T_1}$.
    \item $({}^0\mathcal X_{T_1}, ({}^0\mathcal X_{T_1})^0{})$ is a
    hereditary torsion pair.
    \item There is a perfect Gabriel topology $\mathcal G$ such that
    $\lambda :R\to S$ is equivalent to $\lambda _{\mathcal G}:R\to Q_{\mathcal
    G}$.
\end{enumerate}
\end{PROP}
\Proof The equivalence of $(1)-(3)$ is proved in \cite[2.1]{CTT}.\\
$(2)\Rightarrow (5)$: By Theorem~\ref{teo:hereditarygiraud} we have that the Giraud subcategory $\mathcal X_{T_1}$ is the category $\mathcal X (\mathcal G)$ of $\mathcal G$-closed modules for some Gabriel topology $\mathcal G$. Since
${\mathcal X}_{T_1}$ is a
coreflective subcategory of $\Mod R$ by Lemma~\ref{reflective}, we infer  from
Theorem~\ref{teo:perfectloc}  that $\mathcal G$ is a perfect
Gabriel topology. Then $\lambda$ and $\lambda _{\mathcal G}$ are
in the same epiclass by Theorem~\ref{GP}.\\
$(5)\Rightarrow(4)$: Since $\lambda$ and $\lambda _{\mathcal G}$ are
in the same epiclass,  the perpendicular category $\mathcal X_{T_1}$ and the category $\mathcal X (\mathcal G)$ of all $\mathcal G$-closed modules coincide.
In particular, $\mathcal X_{T_1}$ is a Giraud subcategory, and combining Theorem~\ref{teo:hereditarygiraud}
 and Lemma~\ref{lemma:7}  we know that  the corresponding  hereditary torsion pair is $({}^0\mathcal X_{T_1}, ({}^0\mathcal X_{T_1})^0{})$. \\
$(4)\Rightarrow (3)$: Since $({}^0\mathcal X_{T_1}, ({}^0\mathcal X_{T_1})^0{})$ is a
    hereditary torsion pair,  $^0\mathcal X_{T_1}$ is closed under submodules. Thus $(3)$ is a consequence of  Lemma~\ref{lemma:7}.\fine

\begin{DEF}
We will say that a tilting module \emph{arises from perfect localization} if there is a perfect Gabriel topology $\mathcal G$ such
    that $R$ embeds in $Q_{\mathcal G}$ and $Q_{\mathcal G}\oplus Q_{\mathcal G}/R$ is a tilting module
    equivalent to $T$.
\end{DEF}

\begin{TEO}~\label{teo:tiltingperfectloc}
Let  $T_R$ be a tilting module of projective dimension one. The following
conditions are equivalent.
\begin{enumerate}
    \item There is an exact sequence $\seb{R}{T_0}{T_1}{}{}$ such
    that $T_i\in \Add T$, $\Hom_R(T_1, T_0)=0$ and $\mathcal
    X_{T_1}$ Giraud subcategory of $\Mod R$.
    \item $T$ arises from perfect localization.
\end{enumerate}
\end{TEO}
\Proof
$(1)\Rightarrow (2)$: By Theorem~\ref{teo:tiltinglocalization} there
is an injective ring epimorphism $\lambda: R\to S$ such that $S\oplus
S/R$ is a tilting module equivalent to $T$, and $\mathcal X_{T_1}$ coincides with
the essential image of the restriction functor $\lambda_\ast:\Mod S\to\Mod R$
induced by $\lambda$.
Now, since $\mathcal X_{T_1}$
is a Giraud subcategory,  we infer from   Proposition~\ref{prop:8}  that there exists a
perfect Gabriel topology $\mathcal G$ such that $\lambda _{\mathcal
G}: R\to Q_{\mathcal G}$ and $\lambda$ are in the same epiclass. So,
  $\lambda_{\mathcal G}$ is
injective, and $Q_{\mathcal G} \oplus Q_{\mathcal G}/R$
is a tilting module equivalent to $T$.\\
 $(2)\Rightarrow (1)$: Let $\mathcal G$ be a perfect Gabriel topology  such
    that $R$ embeds in $Q_{\mathcal G}$ and $Q_{\mathcal G}\oplus Q_{\mathcal G}/R$ is a tilting module
    equivalent to $T$. Then the sequence $\seb{R}{Q_{\mathcal G}}{Q_{\mathcal G}/R}{}{}$ has the stated properties. In fact, if $T_1=Q_{\mathcal G}/R$, then we know from Theorem~\ref{teo:tiltinglocalization} that $\lambda _{\mathcal
G}: R\to Q_{\mathcal G}$ induces an equivalence between  $\mathcal X_{T_1}$
and $\Mod Q_{\mathcal G}$. Then $\mathcal X_{T_1}$ coincides with $\mathcal X (\mathcal G)$
 by
Theorem~\ref{teo:perfectloc}, and it is therefore
 a Giraud subcategory.
\fine

 \begin{ES}\label{Fuchs} Exact sequences $\seb{R}{T_0}{T_1}{}{}$ such
    that $T_i\in \Add T$ and $\mathcal
    X_{T_1}$ is  a Giraud subcategory of $\Mod R$ may exist even when $T$ is a tilting module which is not of the form $S\oplus S/R$.

 Let $R$ be a commutative domain, and $Q$ its quotient field.
Denote by $\mathcal D$ the class  of all
{divisible modules}. It was shown by Facchini  that there is a tilting
module of projective dimension one generating $\mathcal D$, namely the {Fuchs' divisible module}
$\delta$, cf.\ \cite[\S VII.1]{FS}.
 Recall further that $\mathcal {D} = \mathcal{U}^{\perp}$ where
${{\mathcal {U}}} = \{ R/rR \mid r \in R \}$ denotes a set of representatives
of  all {cyclically presented modules}.
Moreover, the module $T_1=\delta/R$ in the exact sequence $0\to R\to \delta\to \delta/R\to 0$   is $\Ucal$-filtered, and the perpendicular category $\Xcal_{T_1}=\Xcal_{\Ucal}$ is the class of all divisible torsion-free modules.

Note that the universal localization of $R$ at $\Ucal$ is exactly $Q$, see \cite[3.7]{ASJ}. So  the $\Xcal_{T_1}$-reflection of $R$ is given by the injective flat epimorphism $\lambda:R\to Q$, and $\Xcal_{T_1}$ is a Giraud subcategory of Mod-$R$. On the other hand, $\delta$ has not the form described in Theorem \ref{teo:tiltinglocalization},  unless pd$Q_R\le 1$, that is, $R$ is a Matlis domain, see \cite[2.11 (4)]{ASJ}.

 \end{ES}

 \section{Tilting modules over semihereditary rings}

As we have seen in Remark \ref{finitetype}, the hereditary torsion pair $(\mathcal T,\mathcal F)$ corresponding
to a  perfect Gabriel topology   $\mathcal G$ is always generated by some set of finitely presented modules $\mathcal C$. If the ring $R$ is right coherent, we have a further useful information.

\begin{PROP}\cite[2.8]{H},\cite{K1}\label{He}
Let $R$ be a right coherent ring, and let   $\mathcal G$ be a perfect Gabriel topology on $R$ with associated  hereditary torsion pair $(\mathcal T,\mathcal F)$. Let $\mathcal S$ be the class of all finitely presented modules from $\mathcal T.$ Then   $\mathcal F={\mathcal S}^o$ and  $\mathcal T=\varinjlim \mathcal S.$
\end{PROP}

 We will use this result for comparing perfect localization with universal localization.

 \begin{LEMMA}\label{equalperpe}
Let  $R$ be right coherent, and let $(\mathcal T,\mathcal F)$ be a hereditary torsion pair.  Let $\mathcal S$ be the class of all finitely presented modules from $\mathcal T$, and assume that  $\mathcal T=\varinjlim \mathcal S$. Denote further by $\mathcal C$ the class of all cyclic modules in $\mathcal S$. If  $\mathcal U\subseteq\mathcal S$ satisfies
  $\perpe{\Ucal}\subseteq\mathcal F\cap {\mathcal C}^\perp$, then $\perpe{\Ucal}=\perpe{\mathcal T}$.
  \end{LEMMA}
\Proof
The inclusion
    $"\supseteq"$ follows immediately from the fact that   $\mathcal U\subseteq \mathcal T$. For the reverse inclusion, let $M\in \mathcal X_{\mathcal U}$. Then $M\in \mathcal F$, so we know from   \cite[p. 518-519]{H} that  there is an exact sequence
    $$\seb{M}{E_M}{C_M}{f_{M}}{g_{M}}$$
   where  $E_M\in\mathcal X_{\mathcal T}$ and  $C_M\in \mathcal T$. Then $C_M =\varinjlim S_i$ for some direct system $(S_i)$ in $\mathcal S$.
   We claim that all $\Hom_R(S_i, C_M)=0$. In fact, if $Y$ is a cyclic submodule of $S_i$, then also $Y$ belongs to $\mathcal S$, hence to $\mathcal C$, and therefore  $\Ext{1}{R}{S_i}{M}=0$. This shows that
    every map $h\in \Hom_R(S_i, C_M)$
     factors through $g_M$, thus
    $h=g_{M}h'$ with $h'\in \Hom_R(S_i, E_M)$, and since $E_M\in \mathcal X_{\mathcal T}$, we deduce $h'=h=0$.

    So we conclude  that $C_M=0$ and $M\cong E_M\in \mathcal X_{\mathcal T}$.
    \fine

 Recall that a ring $R$ is said to be  right \emph{semihereditary} if every finitely generated right ideal is projective. Then, by a classical result of Kaplansky, all finitely generated submodules of a  right projective module are projective, hence all finitely presented modules have projective dimension at most one.

     \begin{PROP}\label{perfectisuniversal} Let $R$ be a right semihereditary ring.
Let $\mathcal G$ be a perfect Gabriel topology on $R$, and let  $(\mathcal T,\mathcal F)$
be the hereditary torsion pair associated to $\mathcal G$. Then the ring epimorphism $\lambda_{\mathcal G}:R\to Q_{\mathcal G}$ is equivalent to the universal localization at the set $\Ucal$ of all finitely presented modules from $\mathcal T$.
\end{PROP}
\Proof
 By Proposition \ref{He} we have  $\mathcal F=\mathcal U^o$, and $\mathcal T=\varinjlim \mathcal U$. Denote by $\mathcal C$ the class of all cyclic modules in $\mathcal U$. Of course
$\perpe{\Ucal}\subseteq \mathcal F\cap \mathcal C^\perp$, hence $\perpe{\Ucal}=\perpe{\mathcal T}$
 by Lemma \ref{equalperpe}.

Note that the Giraud subcategory $\mathcal X_{\mathcal T}$ coincides with  the category of $\mathcal G$-closed modules,  see Theorem  \ref{teo:hereditarygiraud}. Thus we infer from Theorem \ref{teo:perfectloc} that
$\mathcal X_{\mathcal T}$ is the essential image of the restriction functor $\Mod Q_{\mathcal G}\to\Mod R$ induced by $\lambda_{\mathcal G}$.
So, it follows from Proposition \ref{teo:univlocalization2} and Theorem \ref{GP} that $\lambda_{\mathcal G}$ and $\lambda _{\mathcal U}$ are in the same epiclass. \fine

\begin{COR}\label{tiltperfectisuniversal} Over a semihereditary ring, every tilting module arising from perfect localization also arises from universal localization.
\end{COR}

 \begin{ES}
  The converse implication  in \ref{perfectisuniversal} or \ref{tiltperfectisuniversal} does not hold true. Indeed,
 \cite[2.2]{CTT} provides an example of a finitely generated tilting module $T$ over a finite dimensional hereditary algebra $R$ admitting
  an exact sequence $\seb{R}{T_0}{T_1}{}{}$ such
    that $T_i\in \Add T$ and  $\Hom_R(T_1, T_0)=0$, but $\mathcal
    X_{T_1}$ is not a Giraud subcategory of $\Mod R$, see also \cite[2.11]{ASJ}.
   Note  that $T$ arises from universal localization at $T_1$, cf.~Example
    \ref{es:pruferartin}.
     \end{ES}

  Let us now focus on the case where $R$ is a Pr\"ufer domain, that is, a commutative semihereditary domain. First of all, we recall the classification of tilting modules
  due to Bazzoni, Eklof and Trlifaj \cite{BET}.

\begin{TEO}{\cite[6.2.15]{GT}}~\label{teo:GT6.2.15}
Let $R$ be a Pr\"ufer domain. There is a bijective correspondence between Gabriel topologies of finite type  and tilting classes.

 The correspondence associates to a Gabriel topology of finite type  $\mathcal L$  the tilting class of all $\mathcal L$-divisible modules. Conversely, if $T$ is a tilting module, then the non-zero finitely generated ideals $I$  such that  $R/I\in  {}^\perp (T^\perp)$ form a basis of the corresponding Gabriel topology.

 \end{TEO}

Over a Pr\"ufer domain, every Gabriel topology of finite type  is perfect.

\begin{LEMMA}~\label{teo:11} Let $R$ be a Pr\"ufer domain. Let further $\mathcal L$ be  a Gabriel topology of finite type, and let $(\mathcal T, \mathcal F)$ be the corresponding hereditary torsion pair. The following statements hold true.

(1) $\mathcal L$ is a   perfect Gabriel topology, $\lambda _{\mathcal L}:R\to Q_{\mathcal L}$ is an  injective ring epimorphism, and $Q_{\mathcal L}/R\in\mathcal T$.

(2) If $Q_{\mathcal L}\oplus
Q_{\mathcal L}/R$ is a tilting module, then the
tilting class $\Gen Q_{\mathcal L}$ coincides with the  class of $\mathcal L$-divisible modules.
\end{LEMMA}
\Proof (1) Either use \cite[XI, 3.5]{St}, or proceed as follows.  By Remark~\ref{finitetype} the torsion class  $\mathcal
T$ consists of all modules $X$ such that every $x\in
X$ has annihilator $\text{ann}_R(x)\in \mathcal L$. In particular  $\text{ann}_R(x)\not=0$, hence  $\mathcal T$
 is contained in the class  of all torsion modules. Thus
$R\in (\mathcal T)^0=\mathcal F$, which
shows that the ${\mathcal X(\mathcal L)}$-reflection $\lambda
_{\mathcal L}:R\to Q_{\mathcal L}$ is injective.

Moreover,   from the construction of the ${\mathcal X(\mathcal L)}$-reflection in \cite[IX, 2.2]{St} or \cite[p. 518-519]{H} we know that
$Q_{\mathcal L}=\{x\in E(R)|\ (R:x)\in \mathcal L\}$ where   $E(R)$
denotes the injective envelope of $R$, and $Q_{\mathcal L}/R\in\mathcal T$. In particular,
$Q_{\mathcal L}$ is an overring of $R$, hence $Q_{\mathcal L}$ is
flat and $\mathcal L$-divisible by \cite[1.1.1, 5.1.15,
5.1.11]{FHP}. So $\mathcal L$ is a perfect Gabriel topology.

(2)
If $Q_{\mathcal L}\oplus
Q_{\mathcal L}/R$ is a tilting module, then $ Q_{\mathcal L}/R$ has projective dimension at most one, and   therefore it has a filtration where the consecutive factors are
 finitely presented cyclic, see  \cite[VI, 6.5]{FS}. Denoting by $\mathcal C$  the class of all cyclic finitely presented modules from $\mathcal T$, we infer that
 $Q_{\mathcal L}/R$ is $\mathcal C$-filtered. Moreover,
 $Q_{\mathcal L}$ is $\mathcal L$-closed and therefore obviously contained in $\mathcal C^\perp$. Then we deduce as in  \cite[3.12]{ASJ}
that
$\Gen Q_{\mathcal L}=\mathcal    C^{\perp}$.

So, it remains to verify that $\mathcal    C^{\perp}$ is the class of $\mathcal L$-divisibles. Now, let $M$ be an $\mathcal L$-divisible module, and let $C=R/I\in\mathcal C$. Then  $I\in \mathcal L$ is  finitely generated, and $M$ is $I$-divisible,  so
we infer from  \cite[6.2.7]{GT} that ${\rm Ext}^1_R(R/I,X)=0$. Conversely, if $M\in\mathcal    C^{\perp}$ and $J\in\mathcal L$, then $J$ contains a finitely generated ideal $I\in \mathcal L$. Since $R/I\in\mathcal C$, we infer again from   \cite[6.2.7]{GT} that $M$ is $I$-divisible, which implies that $M$ is also $J$-divisible.
\fine

Next, we prove a converse of Proposition~\ref{perfectisuniversal}.

\begin{PROP}\label{universalisperfect}
Let $R$ be a Pr\"ufer domain, and let $\Ucal$ be a set of finitely presented cyclic modules. Let further $\mathcal L$ be the Gabriel topology having as basis the set $\mathcal B$ of all non-zero finitely generated ideals $I$  such that  $R/I\in\Ucal$. Then the universal localization at $\Ucal$ is equivalent to $\lambda_{\mathcal L}:R\to Q_{\mathcal L}$.
\end{PROP}
\Proof
The Gabriel topology
$\mathcal L$ is obviously   of finite type, hence perfect by Lemma~\ref{teo:11}. Let  $(\mathcal T,\mathcal F)$
be the hereditary torsion pair associated to $\mathcal L$, and let $\mathcal S$ be the class of all finitely presented modules from $\mathcal T$. By Proposition \ref{He} we have  $\mathcal F=\mathcal S^o$, and $\mathcal T=\varinjlim \mathcal S$.

We verify that $\Ucal$ satisfies the assumptions of  Lemma \ref{equalperpe}.

(i)  $\Ucal\subseteq \mathcal S$. In fact, $\Ucal$ consists of modules of the form $R/I$ where the ideal $I$  belongs to the basis  $\mathcal B$,
so the annihilator of any element of $R/I$  belongs to $\mathcal L$ since it contains $I$. From Remark~\ref{finitetype} we infer   $\Ucal\subseteq \mathcal T$, hence  $\Ucal\subseteq \mathcal S$.

(ii) If $J\in\mathcal L$ and $X\in\perpe{\Ucal}$, then $\Hom_{R}({R/J},{X})=0$. This is because $J$ contains an ideal $I$ such that $R/I\in\Ucal$ and therefore $\Hom_{R}({R/I},{X})=0$.

(iii) $\perpe{\Ucal}\subseteq \mathcal F$. Indeed, if $X\in\perpe{\Ucal}$ and $Y$ is a cyclic module in $\mathcal T$, then $Y=R/J$ with $J\in \mathcal L$, hence $\Hom_R(Y,X)=0$ by (ii). But this implies $\Hom_R(M,X)=0$ for all $M\in\mathcal T$, that is, $X\in\mathcal F$.

(iv) Let $\mathcal C$ be the class of all cyclic modules in $\mathcal S$, and let $X\in\perpe{\Ucal}$. We verify that $X\in\mathcal C^\perp$. If $C\in\mathcal C$, then $C=R/J$ for some finitely generated ideal $J\in\mathcal L$, and $J$ must contain an element from the basis $\mathcal B$, that is, a non-zero finitely generated ideal $I$  such that  $R/I\in\Ucal$. Then ${\rm Ext} ^1_R(R/I,X)=0$, which means
by \cite[6.2.7]{GT}  that $X$ is $I$-divisible. But then $X$ is also $J$-divisible, and again by  \cite[6.2.7]{GT}
we infer ${\rm Ext}^1_R(R/J,X)=0$.

Now Lemma \ref{equalperpe} yields that $\perpe{\Ucal}=\perpe{\mathcal T}$, and we complete the proof as in Proposition
 \ref{perfectisuniversal}.
\fine

 \begin{COR}\label{Giraudsub}
 Let $R$ be a Pr\"ufer domain. Let $T$ be a tilting module of projective dimension one, and let  $\seb{R}{T_0}{T_1}{}{}$ be an   exact sequence where $T_0,T_1\in\Add T$. Then $\perpe{T_1}$ is a Giraud subcategory of $\Mod R$.
 \end{COR}
 \Proof
 By \cite[6.2.10]{GT} there is   a class $\Ucal$ of finitely presented cyclic modules in ${}^\perp (T^\perp)$ such that
 $T_1$ is $\Ucal$-filtered. By Theorem~\ref{teo:univlocalization1} it follows that $\perpe{T_1}$ is the essential image of the restriction functor induced by the universal localization $\lambda_\Ucal$. But $\lambda_\Ucal$ is equivalent to a perfect localization by  Proposition~\ref{universalisperfect}. So, Theorem~\ref{teo:perfectloc} yields that $\perpe{T_1}$
 is a Giraud subcategory.\fine

If $\mathcal L$ is a Gabriel topology of finite type such that the localization $Q_{\mathcal L}$ has projective dimension at most one over $R$, then it was shown by Salce \cite{S2} that the corresponding tilting module $T$ is equivalent to $ Q_{\mathcal L}\oplus Q_{\mathcal L}/R$. We recover Salce's result as a consequence of Theorem~\ref{teo:tiltingperfectloc}. Moreover,  we obtain that every  tilting module of the form $S\oplus S/R$ studied in Theorem~\ref{teo:tiltinglocalization} arises from perfect localization and from universal localization.
\begin{TEO}~\label{teo:14}
Let $R$ be a Pr\"ufer domain. Let $T$ be a tilting module, and let $\mathcal L$ be the  associated Gabriel topology  of finite type. The following statements are equivalent.
\begin{enumerate}
\item $\pdim Q_{\mathcal L}\le 1$.
\item  $T$ arises from perfect localization.
\item $T$ arises from universal localization.
\item There is an  exact sequence $\seb{R}{T_0}{T_1}{}{}$
where $T_0,T_1\in\Add T$ and $\Hom_R(T_1,T_0)=0$.
\end{enumerate}
Moreover, under these conditions,   $T$ is equivalent to $Q_{\mathcal L}\oplus Q_{\mathcal L}/R$.
\end{TEO}
\Proof
First of all, recall that the tilting class $\Gen T$ is the class of all $\mathcal L$-divisible modules.

$(1)\Rightarrow (2)$: We know from Lemma~\ref{teo:11}(1) and Remark~\ref{finitetype}
that $\lambda _{\mathcal L}:R\to Q_{\mathcal L}$ is an  injective ring epimorphism, and that $Q_{\mathcal L}$ is a flat $R$-module. If $\pdim Q_{\mathcal L}\le 1$, then it follows from Proposition~\ref{prop:6} that $Q_{\mathcal L}\oplus
Q_{\mathcal L}/R$ is a tilting module. Since  its tilting class $\Gen Q_{\mathcal L}$ coincides with   the class of $\mathcal L$-divisible modules
 by Lemma~\ref{teo:11}(2), we conclude that $Q_{\mathcal L}\oplus
Q_{\mathcal L}/R$ is equivalent to $T$.

$(2)\Rightarrow (3)$
follows immediately from
Corollary \ref{tiltperfectisuniversal}.

$(3)\Rightarrow (4)$ holds true by Theorem \ref{teo:tiltinglocalization}.

$(4)\Rightarrow (2)$ follows by combining Theorem~\ref{teo:tiltingperfectloc} and Corollary~\ref{Giraudsub}.

$(2)\Rightarrow (1)$:
Let  $\mathcal G$ be a perfect
Gabriel topology such that $\lambda _{\mathcal G}:R\to Q_{\mathcal G}$ is injective and  $Q_{\mathcal G}\oplus Q_{\mathcal G}/ R$ is a tilting module whose tilting class $\Gen Q_{\mathcal G}$ coincides with $\Gen T$.  On the other hand, $\Gen Q_{\mathcal G}$ coincides with the class of $\mathcal G$-divisible modules by Lemma~\ref{teo:11}(2), and Gen$T$   coincides with the class of $\mathcal L$-divisible modules. So, we infer by Theorem~\ref{teo:GT6.2.15} that the Gabriel topologies $\mathcal G$ and $\mathcal L$ coincide. Hence $T$ is  equivalent to $Q_{\mathcal L}\oplus Q_{\mathcal L}/R$, and
 $\pdim Q_{\mathcal L}\le 1$.
 \fine

\end{document}